\documentclass[12pt]{amsart}

\usepackage{geometry}                
\usepackage{graphicx}
\usepackage{amssymb,enumerate}
\usepackage{amsmath,amsfonts,amsthm}\usepackage{exscale}
\usepackage{latexsym}\usepackage{pict2e}
\usepackage{verbatim}\usepackage[all,poly,curve,knot,arrow]{xy}
\usepackage{fancybox}\usepackage{wasysym}
\usepackage{fancyhdr}

\usepackage{epstopdf}
\theoremstyle{plain}
\newtheorem{thm}{Theorem}[section]
\newtheorem{main}[thm]{Theorem}

\newtheorem{lem}[thm]{Lemma}
\newtheorem{prop}[thm]{Proposition}
\newtheorem{cor}[thm]{Corollary}

\newtheorem{ques}[thm]{Question}

\theoremstyle{definition}
\newtheorem{conv}[thm]{Convention}
\newtheorem{ex}[thm]{Example}
\newtheorem{defn}[thm]{Definition}
\newtheorem{rmk}[thm]{Remark}

\newtheorem{fact}[thm]{Fact}

\newcommand{\be}{\begin{enumerate}}
\newcommand{\ee}{\end{enumerate}}
\newcommand{\bd}{\begin{defn}}
\newcommand{\ed}{\end{defn}}
\newcommand{\bp}{\begin{prop}}
\newcommand{\ep}{\end{prop}}
\newcommand{\bl}{\begin{lem}}
\newcommand{\el}{\end{lem}}
\newcommand{\bc}{\begin{cor}}
\newcommand{\ec}{\end{cor}}

\newcommand{\bt}{\begin{thm}}
\newcommand{\et}{\end{thm}}
\newcommand{\bpf}{\begin{proof}}
\newcommand{\epf}{\end{proof}}
\newcommand{\bex}{\begin{ex}}
\newcommand{\eex}{\end{ex}}
\newcommand{\bft}{\begin{fact}}
\newcommand{\eft}{\end{fact}}
\newcommand{\brk}{\begin{rmk}}
\newcommand{\erk}{\end{rmk}}

\newcommand{\ba}{\begin{align*}}
\newcommand{\ea}{\end{align*}}
\newcommand{\tn}{\textnormal}
\newcommand{\df}{\textnormal{df}}
\newcommand{\bit}{\begin{itemize}}
\newcommand{\eit}{\end{itemize}}

\newcommand{\bcm}{}
\newcommand{\cref}[1]{(\ref{#1})}
\newcommand{\ol}{\overline}\newcommand{\wt}{\widetilde}
\newcommand{\hf}{\hfill}
\newcommand{\ci}{\CIRCLE}

\newcommand{\led}[4]{$\xymatrix@1{{#1\,} \ar[q-1]^-{#2}_-{#3}& {\,#4}}$}
\newcommand{\wed}[4]{$\xymatrix@C=27pt@1{{#1\,} \ar@{->>}[q-1]^-{#2}_-{#3}& {\,#4}}$}

\begin{document}

\title[Inequalities on Bruhat graphs, $R$- and KL polynomials]
{Inequalities on Bruhat graphs, $R$- and Kazhdan-Lusztig polynomials}
\author{Masato Kobayashi}
\thanks{This article will appear in Journal of Combinatorial Series A \textbf{120} (2013) no.2, 470--482.}
\date{\today}                                           
\address{Department of Mathematics\\
Tokyo Institute of Technology,
2-12-1 Oookayama, Tokyo 152-8551, Japan.}
\email{kobayashi@math.titech.ac.jp}
\maketitle
\begin{abstract}

From a combinatorial perspective, we establish three inequalities on coefficients of $R$- and Kazhdan-Lusztig polynomials for crystallographic Coxeter groups:
(1) Nonnegativity of $(q-1)$-coefficients of $R$-polynomials,
(2) a new criterion of rational singularities of Bruhat intervals by sum of quadratic coefficients of $R$-polynomials,
(3) existence of a certain strict inequality (coefficientwise) of Kazhdan-Lusztig polynomials.
Our main idea is to understand Deodhar's inequality in a connection with a sum of $R$-polynomials and edges of Bruhat graphs.

\end{abstract}

\tableofcontents

\section{Introduction}\label{sintro}

\begin{figure}[h!]
\caption{Bruhat graph of boolean poset of rank 3} \label{b3}
\[
\begin{xy}
0;<5mm,0mm>:
,0+(-12,0)*{\ci}="b1"*+++!U{}
,0+(-12,2)*{\ci}="am"*+++!U{}
,0+(-12,5)*{\ci}="cm"*+++!U{}
,(-3,2)+(-12,0)*{\ci}="al1"*++!R{}
,(3,2)+(-12,0)*{\ci}="ar1"*++!L{}
,(-3,5)+(-12,0)*{\ci}="cl1"*++!R{}
,(3,5)+(-12,0)*{\ci}="cr1"*++!L{}
,(0,7)+(-12,0)*{\ci}="t1"*+++!D{}
,\ar@{->}"b1";"am"
,\ar@{->}"b1";"ar1"
,\ar@{->}"b1";"al1"
,\ar@{<-}"cr1";"am"
,\ar@{<-}"cm";"ar1"
,\ar@{->}"al1";"cl1"
,\ar@{->}"ar1";"cr1"
,\ar@{->}"cl1";"t1"
,\ar@{->}"cr1";"t1"
,\ar@{->}"cm";"t1"
,\ar@{->}"al1";"cm"
,\ar@{->}"am";"cl1"
\end{xy}\]
\end{figure}
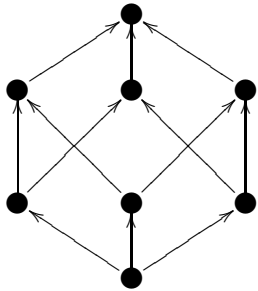

In 1979 \cite{KL}, Kazhdan and Lusztig discovered two families of polynomials (now known as \emph{$R$-} and \emph{Kazhdan-Lusztig polynomials}) in the course of studying Hecke algebras and Schubert varieties. This family of polynomials is indexed by pairs of elements in a Coxeter group, and the polynomials are in one variable and have integer coefficients.
Because Coxeter groups are involved, \emph{Bruhat order} plays a central role in the theory. 
Bruhat order is locally \emph{Eulerian}. Eulerian posets have been of great importance in combinatorics; one particularly important example is that of the face lattices of convex polytopes, and there has been much study of their $f$- and $h$-vectors. We will not list the large number of classical references on this topic but instead refer to books by Stanley \cite{stanley4} and Ziegler \cite{ziegler} and the references therein. \\
\indent Recently, there has been work specifically on the $f$-vectors of lower Bruhat intervals. Bj\"{o}rner-Ekedahl \cite[Theorems A, E]{bjorner3} and Brion \cite[Corollary 2]{brion} have shown certain unimodality properties hold for $f$-vectors of such intervals.
Their approach is of a rather geometric flavor, using the theory of intersection cohomology.\\
\indent From a more combinatorial perspective, the \emph{Bruhat graph}, introduced by Dyer \cite{dyer1}, is one of the most powerful tools for encoding information about Bruhat intervals. Among Bruhat intervals are two classes of fundamental Eulerian structures: \emph{boolean} and \emph{dihedral} intervals.
These coincide up to length 2; however, for length $\ge 3$, their graph structures are different (Figures \ref{ss3}, \ref{b3}). In particular, the graph in Figure \ref{ss3} contains an edge of length 3.
This non-boolean structure leads to the study of labeled Bruhat paths on Bruhat graphs. Dyer \cite{dyer2} gave an interpretation of $\wt{R}$- (and $R$-)polynomials as the generating function of paths with increasing labels in an arbitrary reflection order.
More recently, Billera \cite{billera} and Billera-Brenti \cite{billera2} studied 
Bruhat intervals using quasisymmetric functions that extend the flag $f$- and $h$-numbers. They introduced the complete cd-index as a more sophisticated way to compute $R$- and Kazhdan-Lusztig polynomials.\\
\indent Bruhat graphs, and these polynomials all come into play when we study \emph{rational smoothness} and \emph{singularities} of Bruhat intervals in crystallographic Coxeter groups. Terms ``rationally smooth" and ``singular" come from geometry of Schubert varieties. 
There are many equivalent criteria \cite[Section 13.2]{billey1};
regular Bruhat graphs, trivial Kazhdan-Lusztig polynomials, a boolean-like sum of $R$-polynomials and palindromic Poincar\'{e} polynomials.
Particularly important is \emph{Deodhar's inequality} to which many researchers contributed; Billey \cite{billey2}, Carrell-Peterson \cite{carrell}, Dyer \cite{dyer4}, Kumar \cite{kumar1} and Polo \cite{polo2} in the 1990s.\\
\indent The motivation for this article was to understand Deodhar's inequality in a more explicit connection with a sum of $R$-polynomials: On the one hand, Deodhar's inequality guarantees nonnegativity of a certain integer. On the other hand, $R$-polynomials involve many negative coefficients. The key idea for our approach is to view $R$-polynomials as polynomials in $q-1$, not $q$. Then nonnegativity of $\wt{R}$-polynomials come into the picture as we shall see.
Although this idea is simple, it is useful for analyzing coefficients of not only $R$-polynomials but also Kazhdan-Lusztig polynomials.\\
\indent Our main result consists of three theorems on inequalities of $R$- and Kazhdan-Lusztig polynomials:
\begin{itemize}
\item nonnegativity of $(q-1)$-coefficients of $R$-polynomials (Theorem \ref{nth1}),
\item a new criterion of singularities for Bruhat intervals (Theorem \ref{nth2}),
\item the existence of a strict inequality of Kazhdan-Lusztig polynomials \\(Theorem \ref{nth3}).
\end{itemize}
Proofs are elementary throughout;
Nonetheless, we hope that these results will be some contributions to analysis of such polynomials in the future. 


Here is an organization of the article:\,Sections \ref{snotn} and \ref{srpolys1} record fundamental terminology on Coxeter groups and $R$-polynomials. Section \ref{sr2} gives an explicit description of coefficients for $R$-polynomials with the idea of the absolute length on Bruhat graphs.
In Section \ref{s2df}, we recall a notion of rational smoothness and singularities. In Section \ref{revi}, we give a new interpretaion of Deodhar's inequality in terms of a sum of $R$-polynomials.
After providing a definition and some background on Kazhdan-Lusztig polynomials in Section \ref{skl1}, we prove Theorem \ref{nth3} in Section \ref{skl2}.


\begin{figure}[h!]
\caption{Bruhat graph of dihedral interval of rank 3} \label{ss3}
\[
\begin{xy}
0;<5mm,0mm>:
,0+(-12,0)*{\ci}="b1"*+++!U{}
,(-3,2)+(-12,0)*{\ci}="al1"*++!R{}
,(3,2)+(-12,0)*{\ci}="ar1"*++!L{}
,(-3,5)+(-12,0)*{\ci}="cl1"*++!R{}
,(3,5)+(-12,0)*{\ci}="cr1"*++!L{}
,(0,7)+(-12,0)*{\ci}="t1"*+++!D{}
,\ar@{->}"b1";"al1"
,\ar@{->}"b1";"ar1"
,\ar@{->}"al1";"cl1"
,\ar@{->}"ar1";"cr1"
,\ar@{->}"cl1";"t1"
,\ar@{->}"cr1";"t1"
,\ar@{->}"al1";"cr1"
,\ar@{->}"ar1";"cl1"
,\ar@{->}"b1";"t1"
\end{xy}\]
\end{figure}
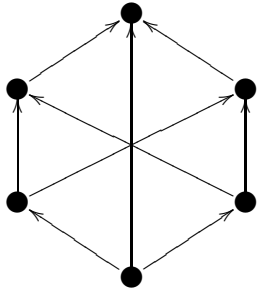

\section{Notation}\label{snotn}
Throughout this article, we follow common notation in the context of Coxeter groups \cite{bjorner2,humphreys}.
By $(W, S)$ (or simply $W$) we mean a Coxeter system with length function $\ell$. Unless otherwise specified, $u, v, w$ are elements of $W$ and $e$ is the unit. Let $T=\cup_{w\in W}\, w^{-1}Sw$ denote the set of reflections. 
Write $u\to w$ if $w=ut$ for some $t\in T$ and $\ell(u)<\ell(w)$. Define \emph{Bruhat order} $u\le w$ if there exist $v_1, \dots, v_n\in W$ such that $u\to v_1\to \cdots \to v_n=w$.
For $u\le w$, let $[u, w]\overset{\tn{def}}{=}\{v\in W\mid u\le v\le w\}$ denote a \emph{Bruhat interval}.
Often $\ell(u, w)\overset{\tn{def}}{=}\ell(w)-\ell(u)$ abbreviates the length of intervals.\\
\indent More notation on polynomials: As usual, the symbol $\mathbb{N}$ indicates the set of nonnegative integers and $\mathbb{Z}$ integers. 
For nonzero $f\in \mathbb{Z}[q]$, say $f$ is \emph{palindromic} if $q^{\deg(f)}f(q^{-1})=f(q)$. Let $[q^n](f)$ denote the coefficient of $q^n$ in $f$. 
An inequality $f\le g $ (or $f\le_q g$) means $[q^n](f)\le [q^n](g)$ for all $n$.
In addition, we use some special notation; see Remark \ref{div}.
\section{$R$-polynomials}\label{srpolys1}

Following \cite[Section 5.1]{bjorner2}, we first give a definition of $R$-polynomials.
\begin{fact}There exists a unique family of polynomials $\{R_{uw}(q)\mid u, w \in W \} \subseteq \mathbb{Z}[q]$ (\emph{$R$-polynomials}) such that
\begin{enumerate}
\item $R_{uw}(q)=0$ \mbox{if $u \not \le w$},
\item $R_{uw}(q)=1$ \mbox{if $u=w$},
\item if $s\in S$ and $ws<w$, then
\begin{align*} R_{uw}(q)=
\begin{cases}
R_{us, ws}(q) &\mbox{ if } us<u,\\
qR_{us, ws}(q)+(q-1)R_{u, ws}(q) &\mbox{ if } u<us.
\end{cases}
\end{align*}
\end{enumerate}
\end{fact}

We can equivalently construct such polynomials from the Hecke algebra of $W$ as in \cite[Chapter 7]{humphreys}. But this definition is enough for our purpose.\\
\indent We will use the following properties of $R$-polynomials later.

\bft{Let $u\le w$.\label{ff}
\be{\item $R_{uw}(q)$ is a monic polynomial of degree $\ell(u, w)$.
\item \label{ntr}If $u\ne w$, then $q-1$ divides $R_{uw}(q)$, i.e., $R_{uw}(1)=0$.
\item $R_{uw}(q)=R_{u^{-1}, w^{-1}}(q)$.
\item We have\label{ralt}
\[ \sum_{ u\le v \le w}(-1)^{\ell(u, v)}R_{uv}(q)R_{vw}(q)=\delta_{uw}\mbox{ (Kronecker delta)}.\]
\item $q^{\ell(u, w)}R_{uw}(q^{-1})=(-1)^{\ell(u, w)}R_{uw}(q)$.\label{fff}
}\ee
}\eft

$R$-polynomials involve many negative $q$-coefficients; However, once we regard them as $(q-1)$-polynomials, we can show the nonnegativity of such coefficients (Theorem \ref{nth1}). 

Next, following \cite[Section 5.3]{bjorner2}, we introduce another family of polynomials associated to $R$-polynomials. They have nonnegative integer coefficients:

\begin{fact}
There exists a unique family of polynomials $\{\widetilde{R}_{uw}(q)\mid u, w \in W \} \subseteq \mathbb{N}[q]$ (\emph{$\widetilde{R}$-polynomials}) such that
\begin{enumerate}
\item $\wt{R}_{uw}(q)=0$ \mbox{if $u \not \le w$},
\item $\wt{R}_{uw}(q)=1$ \mbox{if $u=w$},
\item if $s\in S$ and $ws<w$, then
\begin{align*} \widetilde{R}_{uw}(q)=
\begin{cases}
\widetilde{R}_{us, ws}(q) &\mbox{ if } us<u,\\
\widetilde{R}_{us, ws}(q)+q\widetilde{R}_{u, ws}(q) &\mbox{ if } u<us,
\end{cases}
\end{align*}
\item $\wt{R}_{uw}(q)$ ($u \le w$) is a monic polynomial of degree $\ell(u, w)$,\vspace{.05in}
\item $R_{uw}(q)=q^{\frac{\ell(u, w)}{2}}\wt{R}_{uw}(q^{\frac{1}{2}}-q^{-\frac{1}{2}})$.
\end{enumerate}
\end{fact}

\section{Some nonnegativity of $R$-polynomials}
\label{sr2}

Now the main discussion begins with Bruhat graphs, our central idea.
Recall that $u\to w$ means $w=ut$ for some $t\in T$ and $\ell(u)<\ell(w)$.
\bd{The \emph{Bruhat graph} of $W$ is a directed graph for vertices $w\in W$ and for edges $u\to w$. We can also consider induced subgraphs for subsets of $W$.
By a \emph{Bruhat path} we always mean a directed path (hence a strict increasing chain) $u\to v_1 \to \dots \to v_n=w$ in the Bruhat graph of $W$.
}\ed

\bd{Let $u\le w$. Define the \emph{absolute length} between $u$ and $w$ to be 
\[a(u, w)=\min\{n\ge 0 \mid u\to v_1 \to \dots \to v_n=w\}.\]
}\ed

\brk{Hence 
$u\to w$ is equivalent to $a(u, w)=1$. 
Note that we have $a(u, w)\le \ell(u, w)$ by Chain Property \cite[Theorem 2.2.6]{bjorner2} and furthermore $(-1)^{a(u, w)}=(-1)^{\ell(u, w)}$ since 
$\ell(v_i, v_{i+1})$ is odd at each edge $v_i\to v_{i+1}$.}\erk

\bft\cite[Exercise 35, Chapter 5]{bjorner2}{\label{d1} For all $u$ and $w$, we have
\[
R'_{uw}(1)=\begin{cases}
1 & \mbox{if } u\to w,\\
0 &\mbox{otherwise.}
\end{cases}
\]
Here, $R'_{uw}(1)$ means $\frac{d}{dq}R_{uw}(q)\Bigr|_{q=1}$.
}\eft

\brk{\label{div}In the context of $R$-polynomials, we usually think them as polynomials of integer coefficients. However, it is also possible to regard them as real polynomials so that we can speak of their derivative.
This idea is helpful particularly when we want to compute some specific coefficients:\,
recall from calculus that for given $f(q)\in \mathbb{R}[q]$, $c\in \mathbb{R}$ and an expansion $f(q)= \sum_{n=0}^{d} a_n(q-c)^n$ with $a_n\in \mathbb{R}$, we have $a_n=f^{(n)}(c)/n!$ where $f^{(n)}$ means the $n$-th derivative.
Below, we apply this idea for $R$-polynomials and $c=1$.
For convenience, we adopt special notation: $[(q-1)^n](f)\overset{\tn{def}}{=} f^{(n)}(1)/n!$ for all nonnegative integers $n$. 
In addition, we write $f\le_{q-1} g$ to mean $[(q-1)^n](f)\le [(q-1)^n](g)$ for all $n$.
}\erk

As a consequence of Fact \ref{d1}, we have that if $u\to w$, then $q-1$ divides $R_{uw}(q)$ while $(q-1)^2$ does not. We may ask more: When does $(q-1)^2$ divide $R_{uw}(q)$ in general? What does the rest of $R_{uw}(q)$ other than a power of $q-1$ look like? Below Theorem \ref{nth1} and Corollary \ref{mth1} answer these questions. 
Here we need a lemma:




%


\bl{Let $u<w$,
\label{inc}$a=a(u, w)$ and $\ell=\ell(u, w)$. Then there exist positive integers $c_\ell$, $c_{\ell-2}, \dots, c_a$ such that
\[\wt{R}_{uw}(q)=c_\ell q^\ell+c_{\ell-2}q^{\ell-2}+\dots +c_aq^a.
\]
Consequently, we have
\[R_{uw}(q)=\sum_{k=0}^{\frac{\ell-a}{2}} c_{a+2k}\,q^{\frac{\ell-a-2k}{2}} (q-1)^{a+2k}.\]
}\el

\bpf{For the first statement, refer to \cite[Theorem 2.5]{incitti4}.
Then 
\begin{align*}
R_{uw}(q)&=q^{\frac{\ell}{2}}\wt{R}_{uw}(q^{\frac{1}{2}}-q^{-\frac{1}{2}})\\
&=q^{\frac{\ell}{2}}\sum_{k=0}^{\frac{\ell-a}{2}}c_{a+2k}(q^{\frac{1}{2}}-q^{-\frac{1}{2}})^{a+2k}\\
&=q^{\frac{\ell}{2}}\sum_{k=0}^{\frac{\ell-a}{2}}c_{a+2k}(q^{-\frac{1}{2}}(q-1))^{a+2k}\\
&=\sum_{k=0}^{\frac{\ell-a}{2}}c_{a+2k}\,q^{\frac{\ell-a-2k}{2}}(q-1)^{a+2k}.
\end{align*}
}\epf




\begin{main}
\label{nth1}
Let $u< w$ and $n$ be a nonnegative integer.
If $n<a(u, w)$ or $n>\ell(u, w)$, then $[(q-1)^n](R_{uw})=0$.
Otherwise, $[(q-1)^n](R_{uw})>0$.
In particular, $a(u, w)$ is the largest power of $q-1$ that divides $R_{uw}(q)$. As a consequence, we have $R_{uw}(q)\ge_{q-1}0$ for all $u, w$.
\end{main}

\bpf{Let $a=a(u, w)$ and $\ell=\ell(u, w)$ for simplicity.
Consider the expression of $R_{uw}(q)$ in Lemma \ref{inc}.
Then $q^{\frac{\ell-a-2k}{2}}=\sum\binom{(\ell-a-2k)/2}{i} (q-1)^{i}\ge_{q-1}0$ for all $k$.
As a result, all terms $(q-1)^n$ ($a\le n\le \ell$) appear in the sum with positive coefficients. If $n>\ell$, then $[(q-1)^n](R_{uw})=0$ since $\deg R_{uw}(q)=\ell$.
}\epf

\begin{cor}\label{mth1}
Let $u< w$ and $d=\ell(u, w)-a(u, w)$.
Then there exist unique integers $f_{i-1}$, $h_i$ $(0\le i\le d)$ such that 
\[R_{uw}(q)=(q-1)^{a(u, w)}\left(\sum_{i=0}^df_{i-1}(q-1)^{d-i}\right)=(q-1)^{a(u, w)}\left(\sum_{i=0}^dh_iq^{d-i}\right),
\]
$f_{-1}=h_0=1$ and $f_{i-1}>0, h_i=h_{d-i}$ for all $i$.
\end{cor}

\bpf{Existence of such positive numbers $f_{i-1}$ follows from Theorem \ref{nth1}. 
Since $R_{uw}(q)$ is monic, we have $f_{-1}=h_0=1$.
Observe next that 
\begin{align*}
q^{d}\left(\sum_{i=0}^dh_i(q^{-1})^{d-i}\right)&=\frac{q^{\ell(u, w)}
}{q^{a(u, w)}}\,\frac{R_{uw}(q^{-1})}{(q^{-1}-1)^{a(u, w)}}\\
&=\frac{(-1)^{\ell(u, w)}}{(-1)^{a(u, w)}}\,\frac{R_{uw}(q)}{(q-1)^{a(u, w)}} \tn{ (Fact \ref{ff} (\ref{fff}))}\\
&=\sum_{i=0}^dh_iq^{d-i}.
\end{align*}
The second factor is thus palindromic, i.e., $h_i=h_{d-i}$.
It remains to show that $f_{i-1}$ and $h_i$ are all integers.
For $f_{i-1}$, we can prove by induction on $\ell(w)$: If $\ell(w)=1$, then $u=e, d=0$ so that $f_{-1}=1$. If $\ell(w)\ge 2$, by recursive relations of $R$-polynomials, we may assume that $u<us$ and $ws<w$ for some $s\in S$. Now the inductive hypothesis shows that both $R_{us, ws}(q), R_{u, ws}(q)\ge_{q-1} 0$ with integer coefficients. Therefore so is $R_{uw}(q)$ since
\begin{align*}
R_{uw}(q)&=qR_{us, ws}(q)+(q-1)R_{u, ws}(q)\\
&=(q-1)R_{us, ws}(q)+R_{us, ws}(q)+(q-1)R_{u, ws}(q).
\end{align*}
All $h_i$ are also integers since there exist linear relations 
$h_i=\sum_{j=0}^i(-1)^{i-j}\binom{d-j}{d-i}f_{j-1}$.
}\epf

\brk{Some $h_i$ can be negative (Example \ref{neg}).
We hope to give a combinatorial interpretation of positive integers $f_{i-1}$.
}\erk
Brenti showed the following result \cite[Theorem 6.3]{brenti4}; However, the last statement of Theorem \ref{nth1} now gives a more direct proof.


\bc{Let $u<w$. Then the following are equivalent:\label{c1}
\be{
\item $R_{uw}(q)=(q-1)^{\ell(u,w)}$,
\item $a(u, w)=\ell(u, w)$. In other words, there do not exist $x, y\in [u, w]$ such that $x\to y$ and $\ell(x, y)= 3$.

}\ee
In particular, if $[u, w]$ is boolean, then $R_{uw}(q)=(q-1)^{\ell(u,w)}$.
}\ec

That is, whenever $[u, w]$ contains an edge of length $3$, then $R_{uw}(q)$ has a factor other than $q-1$.
We see a small example.

\bex{\label{neg}
Let $W=\tn{A}_2$, $u=123$ and $w=321$ (one-line notation).
Figure \ref{ss3} shows the Bruhat graph of $[u, w]$. Observe that $u\to w$ with $\ell(u, w)=3$. As computed in \cite[Example 5.1.2]{bjorner2}, the $R$-polynomial of $[u, w]$ is $(q-1)(q^2-q+1)$ . Since $q^2-q+1=(q-1)^2+(q-1)+1$, 
we have 
\[R_{uw}(q)=(q-1)(q^2-q+1)=(q-1)^3+(q-1)^2+(q-1)\ge_{q-1} 0.\]

}\eex

We close this section with one more result; it shows bounds of coefficients of $R$-polynomials by binomial ones.

\bp{
\label{bn}
Let $w\in W$. Then for each $u<w$, we have
\[(q-1)^{\ell}\le_{q-1} R_{uw}(q)\le_{q-1} q^{\ell}.\]
}\ep

\bpf{The first inequality follows from Corollary \ref{mth1}.
For the second, it is enough to show that $[(q-1)^n](R_{uw})\le \binom{\ell}{n}$ for all $n$, $0\le n\le \ell$.
The proof proceeds by induction on $\ell(w)$:
If $\ell(w)=1$, then $R_{uw}(q)=q-1$ so that $[q-1](R_{uw})=1$.
Suppose next $\ell(w)\ge 2$. Choose $s\in S$ such that $ws<w$. If $us<u$, then $R_{uw}(q)=R_{us, ws}(q)$ in which case we are done by induction ($\ell(ws)<\ell(w)$).
If $us>u$, then $R_{uw}(q)=qR_{us, ws}(q)+(q-1)R_{u, ws}(q)$ so that
\begin{align*}
[(q-1)^n](R_{uw})&=[(q-1)^n]((q-1)R_{us, ws}+R_{us, ws}+(q-1)R_{u, ws})\\
&\le \binom{\ell-2}{n-1}+\binom{\ell-2}{n}+\binom{\ell-1}{n-1}\mbox{\quad  (induction)}\\
&= \binom{\ell-1}{n}+\binom{\ell-1}{n-1}\\
&= \binom{\ell}{n}.
\end{align*}

}\epf

\brk{Unfortunately, this is a little different from Brenti's Conjecture: $|[q^n](R_{uw})|\le \binom{\ell}{n}$ \cite[Problem 5.2]{brenti1}. The conjecture remains open at time of writing (March 2012). We hope that our inequality above is helpful for proving it. See also Caselli \cite{caselli2} for some relations between $q$-coefficients of $R$-polynomials and binomial ones.
}\erk


\section{Rational smoothness and singularities}\label{s2df}

In this section, we recall rational smoothness and singularities for Bruhat intervals. This is a key concept in the sequel. We begin with a convention:

\begin{conv}\label{convv}
\emph{In what follows we 
assume that $W$ is crystallographic}, i.e., its Coxeter graph has Coxeter labels only from $\{2, 3, 4, 6, \infty\}$.
\end{conv}

The reason for this assumption is to ensure the correctness of Definition \ref{nons}, Facts \ref{kl2} and \ref{kl22}.




\bd{Let $u\le w$. 
Set \[\ol{N}(u, w)=\{v\in W \mid u\to v\le w\} \mbox{ and } \ol{\ell}(u, w)=|\ol{N}(u, w)|.\]
}\ed

In words, $\ol{N}(u, w)$ is the neighborhood of the bottom vertex on the Bruhat graph of $[u, w]$ (Figure \ref{nuw}); $\ol{\ell}(u, w)$ is the number of those outgoing edges.
\bd{The \emph{defect} of $[u, w]$ is $\df(u, w)=\ol{\ell}(u, w)-\ell(u, w)$.
}\ed 

We know nonnegativity of this integer:
\begin{figure}
\caption{$\ol{N}(u, w)$} \label{nuw}
\[
\begin{xy}
0;<14mm,0mm>:
,(0,4)*+={\ci}*+++!D{w}
,(0,0)*{\ci}*+++!U{u}
,\ar@{->}(0,0);(1,1)
,\ar@{->}(0,0);(-1,1)
,\ar@{->}(0,0);(-0.8,3)
,\ar@{->}(0,0);(0,1)
,\ar@{->}(0,0);(0.8,3)
,\ar@{..}(0,1);(0,4)
,\ar@{..}(-1,1);(0,4)
,\ar@{..}(0,4);(-0.8,3)
,\ar@{..}(0,4);(1,1)
,\ar@{..}(0,4);(0.8,3)
\end{xy}\]
\end{figure}

\bft{{\cite[Deodhar's inequality]{dyer4}}\label{deodd} $\df(u, w)\ge 0$.
}\eft


\bd\cite[Section 13.2]{billey1}{\label{nons} Let $u\le w$. Say $[u, w]$ is \emph{rationally smooth} if we have the following equivalent conditions:
\be{\item $\sum_{x\le v\le w}R_{xv}(q)=q^{\ell(x, w)}$ for all $x$ with $u\le x<w$,\label{nn}
\item $\df(x, w)=0$ for all $x$ with $u\le x<w$.}\ee
Otherwise, say $[u, w]$ is \emph{singular}.\label{sgg}
}\ed
Recall from Theorem \ref{nth1} that $R_{xv}(q)\ge_{q-1} 0$ for all $x, v$.
Hence a sum of such polynomials satisfies the same property.
In this rationally smooth case, we can write the sum in this way:
\[\sum_{x\le v\le w}R_{xv}(q)=q^\ell=(q-1+1)^\ell=\sum_{n=0}^{\ell} \binom{\ell}{n}(q-1)^n.\]
In particular, $[q-1]\left(\sum_{x\le v\le w}R_{xv}\right)=\ell(x, w)$ for $n=1$. In the next section, we establish two results on such coefficients in a more general point of view. They are stated in the same form.

\section{Deodhar's inequality revisited}\label{revi}



%


\bp{\label{dvc}Let $u\le w$. Then for all $x$ with $u\le x<w$, we have
\[[q-1]\left(\sum_{x\le v\le w}R_{xv}\right)-\ell(x, w)\underset{(*)}{\ge} 0.\]
Moreover, $[u, w]$ is singular if and only if \tn{($*$)} is strict for some $x$ with $u\le x<w$.
}\ep
\bpf{Recall from Fact \ref{d1}: $[q-1](R_{xv})=R'_{xv}(1)=\begin{cases}
1&\mbox{if }x\to v,\\ 0 &\mbox{otherwise.}\end{cases}$
It follows that $[q-1]\left(\sum_{x\le v\le w}R_{xv}\right)=|\ol{N}(x, w)|=\ol{\ell}(x, w)$.
Hence ($*$) is nothing but rephrasing Deodhar's inequality: $\df(x, w)=\ol{\ell}(x, w)-\ell(x, w)\ge 0$. 
Consequently, $[u, w]$ is singular if and only if $\df(x, w)>0$ for some $x$ with $u\le x<w$ if and only if $[q-1]\left(\sum_{x\le v\le w}R_{xv}\right)-\ell(x, w)>0$ for some $x$ with $u\le x<w$.
}\epf




\begin{main}\label{nth2}
Let $u\le w$. Then for all $x$ with $u\le x<w$, we have
\[[(q-1)^2]\left(\sum_{x\le v\le w}R_{xv}\right)-
\binom{\ell(x, w)}{2} \underset{(*)}{\ge} 0.\]
Moreover, $[u, w]$ is singular if and only if \tn{($*$)} is strict for some $x$ with $u\le x<w$.


\end{main}

We need three lemmas for the proof of Theorem \ref{nth2}. 

\bl{\cite[Lemma 12.2.12 ($\tn{b}_2$)]{kumar2} If $u\to w$, then we have $(q-1)^2$ divides $R_{uw}(q)-q^{\frac{\ell(u, w)-1}{2}}(q-1)$.\label{le1}}\el

\bl{If $u\to w$, then $c_1=[q](\wt{R}_{uw})=1$ where $c_1$ is the integer as given in Lemma \ref{inc}.\label{le2}
}\el

\bpf{Consider the expression of $R_{uw}(q)$ in Lemma \ref{inc} with $a=1$. Since $(q-1)^2$ divides $R_{uw}(q)-q^{\frac{\ell(u, w)-1}{2}}(q-1)$, we must have $c_1=1$.
}\epf

\bd{By $u\to\to w$ we mean $u<w$ and $a(u, w)=2$.
For such $(u, w)$, define $m(u, w)=|\{v\in [u, w]\mid u\to v\to w\}|$.
}\ed



\bl{\label{le3}\hfill\be{\item If $u\to w$, then $R_{uw}''(1)=\ell(u, w)-1$.
\item If $u\to\to w$, then $R_{uw}''(1)=m(u, w)$. 
}\ee
}\el

\bpf{(1):\,
Suppose $u\to w $.
Consider the expression of $R_{uw}(q)$ in Lemma \ref{inc}.
Differentiate it twice and let $q=1$. Then all terms $k\ge 1$ vanish so that the only $k=0$ term (with $c_1=1$ as above) survives:


\begin{align*}
R_{uw}''(1)&=\left(q^{\frac{\ell(u, w)-1}{2}}(q-1)\right)''\Bigr|_{q=1}=
\frac{4\ell(u, w)-4}{4}=\ell(u, w)-1.
\end{align*}
(2):\,Suppose $u\to \to w$. 
 Differentiate the equation in \mbox{Fact \ref{ff} (\ref{ralt})}
\[
\sum_{u\le v\le w}(-1)^{\ell(u, v)}R_{uv}(q)R_{vw}(q)=\delta_{uw}=0 \label{del}
\]
twice. Then let $q=1$:
\[\sum_{u\le v\le w}(-1)^{\ell(u, v)}(R''_{uv}(1)R_{vw}(1)+2R_{uv}'(1)R_{vw}'(1)+R_{uv}(1)R_{vw}''(1))=0.\]
Note that $R''_{uv}(1)R_{vw}(1)$ is nonzero if and only if $a(u, v)\ge 2$ and $v=w$ (use Fact  \ref{d1} and Corollary \ref{mth1}). 
Similarly $R_{uv}'(1)R_{vw}'(1)$ is nonzero (and must be 1) if and only if $u\to v\to w$. Also $R_{uv}(1)R_{vw}''(1)$ is nonzero if and only if $v=u$ and $a(v, w)\ge 2$. Computing signs, we have $R''_{uw}(1)-2m(u, w)+R''_{uw}(1)=0$. 
}\epf

\bpf[Proof of Theorem \ref{nth2}]{Let $u\le x<w$.
Since $[(q-1)^2](R_{xv}(q))=R_{xv}^{''}(1)/2$ (Remark \ref{div}) for all $v\in [x, w]$, it is enough to show that 
\[\sum_{x\le v\le w}\left(\frac{R_{xv}''(1)}{2}\right)-\binom{\ell(x, w)}{2}\ge 0.\]
In the sum, we only need to consider $v\in [x, w]$ such that $a(x, v)\le 2 \le \ell(x, v)$ (otherwise $R_{xv}''(1)=0$ thanks to Theorem \ref{nth1}). Using Lemma  \ref{le3}, write down the sum separately as
\[\sum_{x\to v\le w}\frac{R_{xv}''(1)}{2}+\sum_{x\to\to y \le w}\frac{R_{xy}''(1)}{2}=\frac{\;1\;}{\;2\;}\left(\sum_{x\to v\le w}(\ell(x, v)-1)
+\sum_{x\to\to y\le w}m(x, y)
\right).
\]
Compute the second term as
\[\sum_{x\to\to y\le w}m(x, y)=\sum_{x\to\to y\le w} |\{v\in [x, y]\mid x\to v\to y\}|
=\sum_{x\to v\le w}\ol{\ell}(v, w).
\]
Now use Deodhar's inequality twice to obtain
\begin{align*}
[(q-1)^2]\left(\sum_{x\le v\le w}R_{xv}\right)&=
\sum_{x\le v\le w}\frac{R_{xv}''(1)}{2}\\
&=\frac{\;1\;}{\;2\;}\sum_{x\to v\le w}\left(\ell(x, v)-1+\ol{\ell}(v, w)\right)\\
&\underset{(**)}{\ge} \frac{\;1\;}{\;2\;}\sum_{x\to v\le w}\left(\ell(x, v)-1+\ell(v, w)\right)\\
&= \frac{\;1\;}{\;2\;}\sum_{x\to v\le w}\left(\ell(x, w)-1\right)\\
&=\frac{\;1\;}{\;2\;}\,\ol{\ell}(x, w)(\ell(x, w)-1)\\
&\underset{(***)}{\ge}\frac{\;1\;}{\;2\;}\,\ell(x, w)(\ell(x, w)-1).
\end{align*}
We thus confirmed the inequality ($*$) in Theorem \ref{nth2} for all $x$ with $u\le x<w$.\\
\indent Suppose moreover that $[u, w]$ is singular. Then 
$\ol{\ell}(x, w)>\ell(x, w)$ for some $x$ with $u\le x< w$ so that ($***$) is strict. Therefore, ($*$) must be also strict.
Suppose, conversely, that ($*$) is strict for some $x$ with $u\le x<w$.
Then ($**$) or ($***$) (or both) must be strict;
equivalently, there exists some $v_0$ such that $x\to v_0\le w$ and $\ol{\ell}(v_0, w)>\ell(v_0, w)$ (hence $v_0\neq w$) or $\ol{\ell}(x, w)>\ell(x, w)$ (or both). Together, we showed that $\ol{\ell}(z, w)>\ell(z, w)$ for some $z$ with $u\le z<w$.
Hence $[u, w]$ is singular.
}\epf

\section{KL polynomials}\label{skl1}

We now turn to Kazhdan-Lusztig polynomials. 
Following \cite[Theorem 5.1.4]{bjorner2}, we first give a definition.
\begin{fact}
\label{kr}
There exists a unique family of polynomials $\{P_{uw}(q)\mid u, w \in W \} \subseteq \mathbb{Z}[q]$ (\emph{Kazhdan-Lusztig polynomials}) such that
\be{
\item $P_{uw}(q)=0$ if $u \not \le w$,
\item $P_{uw}(q)=1$ if $u=w$,
\item $\deg P_{uw}(q)\le (\ell(u, w)-1)/2$ if $u<w$,
\item \label{kr4}if $u\le w$, then
\begin{align*}
q^{\ell(u, w)}P_{uw}(q^{-1})=\sum_{u\le v\le w} R_{uv}(q)P_{vw}(q),
\end{align*}
\item $[q^0](P_{uw})=1$ if $u\le w$,
}\ee
\end{fact}

\bd{\label{sg}Let $u\le w$. Say $u$ (or $[u, w]$) is \emph{singular} if $P_{uw}(q)>1$ where $>$ is the $q$-coefficientwise partial order in $\mathbb{Z}[q]$. Say $u$ is \emph{rationally smooth} if $P_{uw}(q)=1$.
}\ed

\brk{\label{rm1}\hf\be{\item This definition is equivalent to Definition \ref{sgg}; see \cite[Section 13.2]{billey1}.

\item Since $[q^0](P_{uw})=1$ whenever $u\le w$, the condition ``$P_{uw}(q)>1$" is equivalent to $P_{uw}(q)= 1+a_jq^j+\cdots$ for some positive integers $j$ and $a_j$. \label{rm2}}\ee
}\erk








Recall from Convention \ref{convv} that $W$ is crystallographic so that:

\bft[Nonnegativity]{\label{kl2}All coefficients of Kazhdan-Lusztig polynomials in $W$ are nonnegative.
}\eft

\bft[Monotonicity]{\label{kl22}
If $u\le v\le w$ in $W$, then $P_{uw}(q)\ge P_{vw}(q)$; In other words, fixing the second index $w$, the function $P_{-, w}(q)$ on $[e, w]$ is \emph{weakly} monotonically decreasing. 
}\eft 
Historically, these became known first for finite and affine Weyl groups $W$; See \cite[Corollary 4]{irving} for Fact \ref{kl2} and \cite[Corollary 3.7]{braden} for Fact \ref{kl22}.
Further, for example, \cite[Theorem 4.2]{bjorner3} says that these properties hold for all crystallographic $W$.
Then a natural question arises: 
\begin{ques}\label{q1}
Fix $w\in W$. 
For which pair $u<v$ in $[e, w]$, does a strict inequality $P_{uw}(q) > P_{vw}(q)$ occur?\end{ques}

Unfortunately, Fact \ref{kl22} does not tell us anything about this.
The idea is to consider $P_{uw}(1)$: 



%



\bp{\label{mono}
Let $u<v\le w$. Then $P_{uw}(q)> P_{vw}(q) \iff P_{uw}(1)> P_{vw}(1)$.}\ep
\bpf{Suppose $u<v\le w$. Then we have the inequality $P_{uw}(q)\ge P_{vw}(q)$ as assumed above. Say $P_{uw}(q)=1+b_1q+\cdots + b_dq^d$, $P_{vw}(q)=1+a_1q+\cdots +a_dq^d$ with $a_i\le b_i$ for all $i$. If $P_{uw}(q)>P_{vw}(q)$, then $a_j<b_j$ for some $j$ ($1\le j\le d$). Then 
\[P_{uw}(1)-P_{vw}(1)=(b_1-a_1)+\cdots+( b_j-a_j )+ \cdots + (b_d-a_d)>0.\]
In a similar fashion, we can show the converse.}\epf


\brk{In particular, $P_{uw}(1)\ge P_{ww}(1)=1>0$ whenever $u\le w$. These positive integers $\{P_{uw}(1)\}$ play an important role in representation theory of Verma modules. This is one of the reasons we want to study it. Here we refer to only \cite[Chapter 8]{humphreys3} in this direction.}\erk

Now, keeping Proposition \ref{mono} in mind, let us put Question \ref{q1} this way:
\begin{ques}\label{q5}
Fix $w\in W$. Further, let $u$ be an arbitrary but fixed element in $[e, w]$ such that $P_{uw}(1)>1$.
Then, for which $v$ in $[u, w]$, does a strict inequality $P_{uw}(1) > P_{vw}(1)$ occur?\end{ques}

Clearly, this is the case for $v=w$ since $P_{ww}(1)=1$. However, we would like to find some $v$ closer to $u$.
Since Bruhat order is defined as the transitive closure of edge relations, it is meaningful to first consider vertices incident to $u$ in $[u, w]$ (Figure \ref{nuw}).
For convenience, let us introduce the following definition:

\bd{An edge $u\to v$ in $[u, w]$ is 
\emph{strict} if $P_{uw}(1)>P_{vw}(1)$.
}\ed
Now, suppose $P_{uw}(1)>1$. Is $u$ incident to some strict edge?
Theorem \ref{nth3} asserts that this is the case for \emph{every} singular vertex $u$ under $w$.





\section{Existence of a strict inequality of KL polynomials}

\label{skl2}
Before Theorem \ref{nth3}, we need a lemma:

\bl{\label{lm}Let $u\le w$. Then
\be{\item we have 
\[\ell(u, w)P_{uw}(1)-2P'_{uw}(1)
=\sum_{v\in \ol{N}(u, w)}P_{vw}(1).\]
\item if $u$ is singular, then $-2P'_{uw}(1)<0$.
}\ee
}\el

\bpf{(1) Differentiate the equation in Fact \ref{kr} (\ref{kr4}) once and let $q=1$. Then the right hand side is $\sum_{v\in \ol{N}(u, w)}R_{uv}'(1)P_{vw}(1)$ thanks to Fact \ref{d1}.
 (2) follows from nonnegativity of coefficients (Fact \ref{kl2}).
}\epf



\begin{main}
\label{nth3}
Let $u\le w$. 
If $P_{uw}(1)>1$, then there exists $t\in T$ such that 
\[P_{uw}(1)>P_{ut, w}(1)>0.\]
\end{main}

\bpf{
Let $n\overset{\tn{def}}{=}|\{v\in \ol{N}(u, w) \mid u\to v \tn{ is strict}\}|$. Suppose $n\le \df(u, w)$.
Then Lemma \ref{lm} implies that
\begin{align*}
\ell(u, w)P_{uw}(1)-2P'_{uw}(1)&=\sum_{v\in \ol{N}(u, w)}P_{vw}(1)\\
&=\sum_{\substack{{v\in \ol{N}(u, w)}\\\tn{strict}}}P_{vw}(1)+(\ol{\ell}(u, w)-n)P_{uw}(1).
\end{align*}
Thus we have
\begin{align*}
\underbrace{-2P'_{uw}(1)}_{< 0}&=\sum_{\substack{{v\in \ol{N}(u, w)}\\\tn{strict}}} P_{vw}(1)+(\ol{\ell}(u, w)-n-\ell(u, w))P_{uw}(1)\\
&=\underbrace{\sum_{\substack{{v\in \ol{N}(u, w)}\\\tn{strict}}} P_{vw}(1)}_{\ge 0}+\underbrace{(\df(u, w)-n)P_{uw}(1)}_{\ge 0},
\end{align*}
a contradiction. Therefore $n\ge \df(u, w)+1\ge 1$.
}\epf

We can repeat this argument as long as $P_{ut, w}(1)>1$ as in the following observation:
\begin{cor}
From every singular vertex $u$ under $w$, there exists a directed path
\[u=v_0\to v_1\to v_2\to \cdots \to v_d \,\,(\le w)\] 
such that $d\ge 1$, all $v_i\to v_{i+1}$ are strict and $v_d$ is rationally smooth.
\end{cor}

\bpf{Suppose $u$ is singular under $w$. As shown in Theorem \ref{nth3}, there exists a strict edge under $w$, say $u\to v_1$. If $v_1$ is rationally smooth, then we are done. Otherwise find another strict edge, say $v_1\to v_2$.
Continue this algorithm until our directed path arrives at some rationally smooth vertex.
}\epf







\vspace{1ex}
\begin{center}
\textbf{Acknowledgments.}\\
\end{center}
I thank the editor and anonymous referees for many helpful comments and suggestions to improve the manuscript.


\bibliography{Kobayashi_rpolys}
\bibliographystyle{amsplain}

\end{document}